\documentclass[11pt]{crmp-l}
\usepackage{amsmath,amsfonts,amsthm,amssymb,amscd,url}
\def\classification#1{\def\@class{#1}}
\classification{\null}

\DeclareFontFamily{OT1}{rsfs}{}
\DeclareFontShape{OT1}{rsfs}{n}{it}{<-> rsfs10}{}
\DeclareMathAlphabet{\mathscr}{OT1}{rsfs}{n}{it}

\DeclareMathOperator{\mo}{\,mod}

\DeclareMathOperator{\Var}{Var}

\DeclareMathOperator{\Prob}{Prob}
\DeclareMathOperator{\Gal}{Gal}

\DeclareMathOperator{\entropy}{entropy}

\newtheorem*{main}{Main Theorem}

\title[Power-free values]{Power-free values, repulsion between points, differing beliefs and the
  existence of error}
\author{Harald Andr\'es Helfgott}
\address{H. A. Helfgott, Mathematics Department, University of Bristol,
Bristol, BS8 1TW, United Kingdom}
\email{h.andres.helfgott@bristol.ac.uk}
\begin{document}
\begin{abstract}
Let $f$ be a cubic polynomial. Then there are infinitely many primes $p$
such that $f(p)$ is square-free.
\end{abstract}
\maketitle
An integer $n$ is said to be {\em square-free} if it is not divisible by
any squares other than $1$. More generally, $n$ is {\em free of $k$th powers}
if $d\in \mathbb{Z}, d^k|n \Rightarrow d = \pm 1$; square-freeness is what
we get in the case $k=2$. Being square-free -- or
at least free of $k$th powers for some $k$ -- 
is a desirable property: many things are
easier to prove for square-free numbers. Thus, given a set of integers,
it is good to know whether infinitely many of them -- or a positive proportion
of them -- are square-free\footnote{If a technique is strong enough to 
prove that infinitely many numbers in the bag are square-free, 
it is generally also strong enough to show that a positive proportion are,
and even to show which proportion are divisible by certain specific squares
and no others. This is certainly the case for all the techniques discussed
here. Results of this strength are necessary
for applications in which, for example, we need to go from the 
discriminant of an 
elliptic curve to its conductor, which is essentially the product of the
prime factors of the discriminant. See \cite{Hsq} for some general
machinery.}.

Let $f$ be a polynomial with integer coefficients. Is there an infinite
number of integers $n$ such that $f(n)$ is free of $k$th powers? 
There are some polynomials for which the answer is clearly ``no'': say
$k=2$ and $f(n) = 4 n$ or $f(n) = n^2$. Slightly more subtly, consider
$f(n) = n (n+1) (n+2) (n+3) + 4$, which is always divisible by $4$. 
Assume, then, that $f$ has no factors repeated $k$ times and that
the following local condition holds: 
\begin{quotation}
(*)\;\;\; for every prime $p$, $f(x) \not\equiv 0 \mo p^k$ has at least one solution in 
$\mathbb{Z}/p^k \mathbb{Z}$.
\end{quotation}
 (Both conditions are obviously necessary, and
both can be checked easily in bounded time.) Then, it is believed,
there must be an infinite number of integers $n$ such that $f(n)$ is
free of $k$th powers.

If $\deg(f) \leq k$, it is not hard to prove as much. If $\deg(f)>k+1$,
the problem is too hard, at least when $k$ is small. For $\deg(f)=k+1$,
the statement was proven by Erd\H{o}s in 1953. In particular, if
$f$ is a cubic polynomial without repeated factors and $f$ satisfies
the local condition (*), then there are infinitely many integers $n$ such
that $f(n)$ is square-free.

Like many proofs in analytic number theory, Erd\H{o}s's proof is rather
tricky, in that it uses the fact that most integers are not prime in order
to avoid certain essential issues. Perhaps because of this, Erd\H{o}s asked:
for $f$ cubic, 
are there infinitely many primes $p$ such that $f(p)$ is square-free?
(More generally, for $f$ of degree $k+1$, are there infinitely many
primes $p$ such that $f(p)$ is free of $k$th powers?) He conjectured that
there are, but, as might be expected, most tricks used before break down.

Hooley (\cite{Ho}, \cite{Ho2}) proved Erd\H{o}s's conjecture
for $k\geq 51$. At about the same time, Nair \cite{Na} proved
the conjecture for $k\geq 7$, using rather different methods. (He was also
the first to treat polynomials with $\deg(f) \geq k+2$, $k$ large;
Heath-Brown (\cite{HB}) has attained further progress in this front.) 
In both approaches, $k$ small is harder than $k$ large; in particular,
the case $k\leq 6$ remained open.

In \cite{He}, I proved Erd\H{o}s's conjecture for all polynomials $f$
with {\em high entropy}; in particular, the proof works when $k=2$,
$\deg(f)=3$ and $\Gal(f) = A_3$. However, most cubic polynomials have
Galois group $S_3$, and their case remained open until now.

I have now managed to prove Erd\H{o}s's conjecture for general cubics.
\begin{main}
Let $f\in \mathbb{Z}\lbrack x\rbrack$ be a cubic polynomial without repeated
roots. Assume that, for every prime $q$, there is a solution 
$x\in (\mathbb{Z}/q^2 \mathbb{Z})^*$ to $f(x) \not\equiv 0 \mo q^2$.
Then there are infinitely many primes $p$ such that $f(p)$ is square-free.
\end{main}
In fact, $f(p)$ is square-free for a positive proportion $C_f$ of all primes $p$
-- where $C_f$ is exactly what we would expect (an infinite
 product of local densities).

The tools used, developed and sharpened in the proof are mostly from
diophantine geometry and probabilistic number theory; there is a key use of 
the modularity of elliptic curves. Let us take a quick walk through the proof.
(A full account shall appear elsewhere.)

It is not hard to show that
\begin{equation}\label{eq:ogto}\begin{aligned}
|\{p\leq N:\; &\text{$f(p)$ is square-free}\}| = C_f \frac{N}{\
\log N} \cdot (1 + o(1)) \\ &+ O(|\{x,y\leq N, d\leq N (\log N)^{\epsilon} : \text{$x$, $y$ prime},\;
d y^2 = f(x)\}|),\end{aligned}\end{equation}
where $|S|$ is the number of elements of a set $S$ and
$C_f$ is the (non-zero) constant we would expect. This (well-known) initial
step goes roughly as follows.
Small square factors
can be sieved out because we know how many primes there are in arithmetic
progressions to small moduli; medium-sized square factors amount to a small
error term, since $\sum_{d>m} N/d^2$ is quite small with respect to $N$
 as soon as $m$ is
moderately large. Large square factors cannot be brushed aside by the same
argument as medium-sized square factors 
simply because there are so many of them:
an additional term that is overshadowed by $N/d^2$ in the medium range
comes to the fore here. This is why the contribution of the
large square factors figures in (\ref{eq:ogto}) as the error
term within $O(\cdot)$.
The sole problem from now on, then, is to show that the expression within $O(\cdot)$
is $o(N/\log N)$.

As you can tell from the notation, we intend to see this as a problem of
bounding the number of integer points $(x,y)$ on curves $d y^2 = f(x)$,
$f$ a fixed cubic polynomial. The issues are two. First, we need very good bounds -- almost
as good as $O(1)$ for the number of points for each $d$, or at least for
every typical $d$. Second, even a bound of $O(1)$ would not be enough!
There are $N (\log N)^{\epsilon}$ curves to consider, and a bound of $O(1)$
per curve would amount to a total bound of $O(N (\log N)^{\epsilon})$, whereas
we need $o(N/\log N)$.

Let us begin with the first issue: we want good bounds on the number of 
integer points $(x,y)$, $x,y\leq N$, on the curve $C_d$ described by
$d y^2 = f(x)$, $d$ fixed. Most techniques for bounding the number of
points on curves are based on some kind of {\em repulsion}: if there
are bees in a room, and each bee stays a yard or more away from every other 
bee, there cannot be too many bees in the room. Repulsion may happen in
the visible geometry of the curve, viz., its graph, as in \cite{BP}; such
a perspective, unfortunately, would not be nearly enough in our case. 
Alternatively, we may look at repulsion in the Mordell-Weil lattice 
corresponding to the curve.

{\em Rational and integer points on curves.}
Let $C$ be a curve of genus $g>0$ 
over $\mathbb{Q}$ (or a number field $K$). The curve
$C$ can be embedded in its Jacobian $J_C$. The rational points
$J_C(K)$ in the Jacobian form a finitely generated abelian group
under the group law of the Jacobian; they are, furthermore, endowed with
a natural norm given by the square root of the canonical height. Hence 
$J_C(K)$ can be naturally embedded in $\mathbb{R}^r$, where
$r$ is the rank of $J_C(K)$. We thus have an (almost) injective map
\[\iota:C(K) \to L \subset \mathbb{R}^r,\]
where $C(K)$ is the set of rational points on $C$ and $L$
is a lattice in $\mathbb{R}^r$. What can be said about 
the image $\iota(C(K))$?

If the genus $g$ is $1$, $\iota(C(K))$ is all of $L$. However,
if $g>1$, then $\iota(C(K))$ looks quite sparse within $L$. 
Mumford \cite{Mu} proved that the points of $\iota(C(K))$ repel each
other: for any $P_1,P_2\in \iota(C(K))$ at about the same distance
from the origin $O$, the angle $\angle P_1 O P_2$ separating $P_1$ from $P_2$
is at least $60^{\circ}$ (for $g=2$), $70.5^{\circ}$ (for $g=3$),
$75.5^{\circ}$ (for $g=4$), \dots -- in general, 
$\angle P_1 O P_2 \geq \arccos \frac{1}{g}$. 

Assume now that the points $P_1$ and $P_2$ come from {\em integer} points
on $C$. Then, as I pointed out in my thesis (\cite[Ch.\ 4]{Heth} or
\cite[Lem.\ 4.16]{Hesq}; see also the earlier work of Silverman
\cite{Si}, \cite[Prop.\ 5]{GS}, which seems to have originated from
the same observation) the angle $\angle P_1 O P_2$ is larger
than if $P_1$ and $P_2$ were merely rational: the angle is at least
$60^{\circ}$ for $g=1$, $75.5^{\circ}$ for $g=2$, \dots -- in general, at least
$\arccos \frac{1}{2 g}$. We obtain better bounds as a consequence. (If
$g=1$, we obtain bounds where Mumford's work does not by itself give any.) 
The bounds are obtained by means of results on sphere-packing; indeed, the number of points fitting at a certain distance from 
the origin and
at a separation of $\geq 60^{\circ}$ from each other is precisely the
number of solid spheres that can fit around a sphere of the same size
in the given dimension.

In the case of the curve $E_d:d y^2 = f(x)$, the bounds we obtain are of
the form $c_1^{\omega(d)}$, $c_1>1$ fairly small ($<2$). This is still not good 
enough, as, on the average, 
it amounts to about $(\log N)^{c_2}$, $c_2$ a small 
but fixed non-zero constant,
for $d\sim N$; what we would like is a bound of the form $(\log N)^{\epsilon}$.

{\em Visible vs.\ Mordell-Weil geometry.} In \cite{HV}, we remark upon the
following phenomenon. Let $P_1$, $P_2$ be two integer (or rational) points
on $C(K)$ at about the same distance from the origin. Suppose that
their coordinates $(x_1,y_1)$, $(x_2,y_2)$ are close to each other, either
in the real place (that is, $|x_1-y_1|$ and $|x_2-y_2|$ are small) or
$p$-adically (that is, $x_1 \equiv x_2 \mo d$ and $y_1 \equiv y_2 \mo d$
for a large integer $d$). Then the angle $\angle P_1 O P_2$ in the
Mordell-Weil lattice is even larger
than it would already have to be. In other words: if two points are close
to each other in the graph of the curve in $\mathbb{R}^2$, they must
be especially far from each other in the Mordell-Weil lattice. 
Thus, if
we partition the set of all rational points into sets of points close
to each other in the graph of the curve, we shall obtain an especially
good bound on the number of elements of each such set. The main concern
is then to keep the number of such sets small.

In the case of the curve $E_d : d y^2 = f(x)$, we have that any two
points $(x_1,y_1)$, $(x_2,y_2)$ on $E_d$ 
induce points $P_1 = (x_1, d^{1/2} y_1)$,
$P_2 = (x_2, d^{1/2} y_2)$ on $E : y^2 = f(x)$. The $y$-coordinates of
$P_1$ and $P_2$ are already close to each other modulo $d$ (that is, 
modulo the prime ideals in $\mathbb{Q}(d^{1/2})$ dividing $d$). The 
congruence classes mod $d$ into which $x_1$ and $x_2$ may fall are
rather restricted, as $f(x_1) \equiv 0 \mo d$ and $f(x_2) \equiv 0 \mo d$;
the total number of congruence classes $x$ modulo $d$ for which
$f(x) \equiv 0 \mo d$ is at most $3^{\omega(d)}$. If $P_1$ and $P_2$
have $x$-coordinates in the same congruence class modulo $d$, then the
angle $\angle P_1 O P_2$ turns out to be at least $90^{\circ}$, or
$90^{\circ} - \epsilon$. Very few points can fit in $\mathbb{R}^r$ lying
at about the same distance from the origin and subtending angles of
$90^{\circ}-\epsilon$ or more from each other. 

There is the problem that $3^{\omega(d)}$ is too large -- a power of
$(\log N)$, on the average, since $\omega(d)$ is usually about $\log \log N$.
However, with probability $1$, an integer $d\leq N$ has a large divisor $d_0$ 
($>N^{1 -
  \epsilon}$) with few prime divisors ($< \epsilon \log \log N$). We may
thus consider points $P_1$, $P_2$ congruent to each other modulo $d_0$,
rather than $d$, and obtain angles $\angle P_1 O P_2$ of size at least
$90^{\circ}-\epsilon$ while considering at most $3^{\omega(d_0)}$ possible 
congruence classes. The total bound on the number of integer points
$(x,y)$ on $E_d:d y^2 = f(x)$ with $x,y\leq N$ is $O((\log N)^{\epsilon'})$
for every typical $d$, that is, for each $d\leq N$ outside a set
of size at most $N/(\log N)^{1000}$.

This is almost as good as $O(1)$. The problem, as said before, is that
this is not good enough; since there are $N$ integers $d=1,2,\dotsc,N$ to consider,
the total bound would be $O(N)$. The issue, then, is how to eliminate most
$d$. Probabilistic number theory has just made its first appearance; 
it shall play a crucial role in what follows.
  
{\em The perspective of the value and the perspective of the argument.}
Our task is still to show that
\begin{equation}\label{eq:malato}
|\{x,y\leq N, d\leq N (\log N)^{\epsilon} : \text{$x$, $y$ prime},\;
d y^2 = f(x)\}|\end{equation}
is at most $o(N/\log N)$. We have a good bound for each $d$, namely,
$O((\log N)^{\epsilon})$ for each $d$ outside a small set,
and a reasonable bound for each $d$ inside that small set.
The idea will now be to consider, in a solution to $d y^2 = f(x)$, what
kind of integer $d = f(x)/y^2$ typically is, and whether it looks much like
a typical integer $d$. We will show that, for most $x$, the integer
$d = f(x)/y^2$ must look rather strange, and that thus there can be
few such $d$.
Stated otherwise: we shall prove that
every prime $x\leq N$ must either lie within a fixed set of size $o(N/\log N)$
 or be such that, if $d y^2 = f(x)$ for some prime $y$ and some integer 
$d\leq N (\log N)^{\epsilon}$, then $d$ must lie within a fixed set of size
$O(N/(\log N)^{1 + 10 \epsilon})$. Combined with our bound for each $d$,
this will yield immediately that (\ref{eq:malato}) is indeed
$o(N/\log N)$.
  
What are, then, the ways in which $f(x)$ will tend to be strange for
a random prime $x$? And which of those ways of strangeness will carry
over to $d$, if $f(x)$ can be written in the form $d q^2$, where $q$
is a large prime?

As far as the second question is concerned: since $q$ is prime,
$f(x)$ and $d$ have almost
the same number of prime divisors. Thus, if we can show that
the number of prime divisors of $f(x)$ is strange for $x$ random, we will
have shown that the number of prime divisors of $d$ is strange for
$x$ random.

Now, for $x$ random, the number of prime divisors $w(f(x))$ will be
about $\log \log N$. Thus, $w(d)$ will also be about $\log \log N$.
Unfortunately, this is typical, not strange, for an integer of the size of
$d$.

Consider, however, primes of different kinds.
Let $K = \mathbb{Q}(\alpha)$, where $\alpha$ is a root of $f(x) = 0$.
Then some primes $p$ will split completely in $K/\mathbb{Q}$,
some primes will not split at all, and some primes may split yet not
split completely. We can write $w_1(n)$, $w_2(n)$ and $w_3(n)$ for
the number of prime divisors of $n$ of each kind. Then, as we shall
see, $w_j(f(x))$ (and thus $w_j(d)$) will tend to be strange for
$x$ random.

Suppose first that $K/\mathbb{Q}$ has Galois group $A_3$. Then
every prime $p$ must either split completely or not split at all.
If $p$ does not split at all, then $f(x) \equiv 0 \mo p$ has no
solutions. Hence $w_2(f(x)) = 0$, and so $w_2(d)=0$. This is certainly
atypical for an integer $d\leq N$. (Usually $w_2(d) \sim \frac{2}{3} 
\log \log N$.) Now suppose $p$ splits completely. Then
$f(x) \equiv 0 \mo p$ has three solutions mod $p$. Thus, for a random
prime $x$, we shall have $f(x) \equiv 0 \mo p$ with probability $3/p$.
Hence $w_2(f(x))$ will most likely be about $\sum_{\text{$p$ splits
completely}} \frac{3}{p} \sim \log \log N$. Thus
$w_2(d) \sim \log \log N$, whereas an integer $d\leq N$ usually
has $w_2(d) \sim \frac{1}{3} \log \log N$.

It is not enough, however, to show that $d$ is strange (i.e., in a set
of size $o(N)$); we must show that $d$ is strange enough (i.e., in a set
of size $O(N/(\log N)^{1 + \epsilon})$). How odd is it for a random integer
$d\leq N$ to have $w_1(d)=0$ and $w_2(d)\sim \log \log N$? The
number
$w_1(d)$ equals $\sum_{\text{$p$ does not split}} X_p$, 
where $X_p$ is a random variable
taking the value $1$ when $p|d$ and $0$ otherwise. 
Similarly, $w_2(d) = \sum_{\text{$p$ splits completely}} X_p$.
Now $X_p$ is $1$ with probability $1/p$ and $0$ with probability $1-1/p$.
Suppose the variables $X_p$ to be mutually independent. Then the
theory of large deviations (Cramer's theorem, or, more appropriately,
Sanoff's theorem) offers the upper bound
\[\Prob\left(\sum_{\text{$p$ does not split}} X_p = 0\; \wedge
  \sum_{\text{$p$ splits completely}} X_p > (1- \epsilon) \log \log N\right)
\ll \frac{1}{(\log N)^{\log 3 - \epsilon'}}.\]

Now, of course, the variables $X_p$ are not actually mutually independent;
the variables $X_{p_1},X_{p_2},\dotsc,X_{p_k}$ can be assumed to be
(approximately) mutually independent only when $p_1 p_2 \dotsb p_k < N$.
However, the fact that $X_p$ has only a small probability of being non-zero
allows us to use the main
 technique from Erd\"os and Kac's Gaussian paper \cite{EK} 
to show that we may treat the variables
$X_p$, for our purposes, as if they were mutually independent. Thus we do
obtain
\[\Prob\left( w_1(d) = 0\; \wedge w_2(d) > (1 - \epsilon) \log \log N\right)
 \frac{1}{(\log N)^{\log 3 - \epsilon'}} = O(N/(\log N)^{1 + \epsilon''}),\]
as desired.

We are done proving the main theorem when $\Gal_f = A_3$. What happens
when $\Gal_f = S_3$? While our analysis is in the main still valid,
the exponent that we obtain instead of $\log 3$ is $\frac{1}{2} \log 3$,
which is less than $1$, and thus insufficient. (In general, for
$d y^k = f(x)$, $\deg(f) = k+1$, the exponent we get may be expressed
as an {\em entropy}, which will depend on $\Gal(f)$ alone. Sometimes
$\entropy(\Gal(f))>1$, and we are done, and sometimes, as in the case of 
$\Gal(f)=S_3$, the entropy is $<1$.)

The reason is that, when $\Gal_f = S_3$, half of the primes split in
$K/\mathbb{Q}$. These primes divide $f(x)$ with exactly the same probability
that they would divide a random integer, and thus they are useless.
What is to be done, then? How can one bridge a gap of size
$1/(\log N)^{1 - \frac{1}{2} \log 3}$?

{\em The existence of error. Modularity.}
Again: what is a way of strangeness such that, if $f(x)$ is strange and
$d = f(x)/q^2$, $q$ a prime, then $d$ must be strange as well? Having
too few or too many prime factors of some kind is one way. Is there
another one?

We have used the fact that $q^2$ has only one prime factor; let us now
use the fact that $q^2$ is a square. For any prime modulus $p$, the
integer $d$ will be a square mod $p$ if and only if $f(x)$ is a square
mod $p$. Now, a random integer is as likely to be a square mod $p$ 
as a non-square
mod $p$. How likely is $f(x)$ to be a square mod $p$ for a random integer $x$
(or a random prime $x$)?

By the Weil bounds, there are $p + O(\sqrt{p})$ points on the curve 
$y^2 = f(x) \mod p$. Hence the probability that $f(x)$ will be a square
mod $p$
is $\frac{1}{2} + O(p^{-1/2})$. This is not good, as $\frac{1}{2}$
would be the probability if there were nothing amiss to be exploited.
Let us show that an error of size about 
$p^{-1/2}$ is in fact present a positive proportion of the time.

Write the number of points on the curve
$y^2 = f(x) \mod p$ as $p + 1 - a_p$. Then the probability that
$f(x)$ will be a square $\mo p$ is precisely
$\frac{1}{2} - \frac{a_p}{2 p} + O(1/p)$; we have to give a lower bound,
on the average, for the size of $|a_p|/2p$ (or, rather,
$a_p^2/p^2$, since we shall later use a variance bound). 
Now, the $L$-function of $E:y^2 = f(x)$ is
$\sum a_n n^{-s}$. By the modularity of elliptic curves (proven by Wiles
et al.), there is a modular form $g$ associated to $L$. We may, in turn,
define a Rankin-Selberg $L$-function $L_{g\otimes g}$
associated to $g$, and use the standard facts that
$L_{g\otimes g} = \sum a_n^2 n^{-s}$ and that $L_{g\otimes g}$ has a simple
pole at
$s=2$. By some Tauberian work (or proceeding as in the proof of
 the prime number theorem)
we deduce that
$\sum_{p\leq z} |a_p|^2/p^2$ is asymptotic to $\log \log z$; in other words,
$a_p^2$ is of size about $p$ on the average.

It is somewhat unpleasant to have to use modularity here, as we need not
know the behaviour of $L_E$ (or $L_{g\otimes g}$) inside the critical strip.
Still, it is hard to see how to do without modularity or some strong kindred
result.
Marc Hindry and Mladen Dimitrov have pointed out to me that, if one wants
to give a (conditional) statement on $k$th-power-free values of polynomials
of degree $k+1$, $k>2$, it may be simpler and more proper to work assuming
Tate's conjecture on $L_{C\times C}$ rather than automorphicity.

{\em Using many small differences. Exponential moments and high moments.}
Now, how may we use these small differences between the probability of
$d$ being a square (for $d$ a random integer) and the probability of
$d$ being a square (for $d = f(x)/q^2$, $x$ a random prime)?

Suppose I am throwing a fair coin in the air. A gentle wind blows;
it may change directions very often, but becomes gradually milder.
I know that the wind has a slight effect on the way the coin lands: if the
wind blows from the east, then, I posit, heads are more likely, whereas, if
the wind blow from the west, tails are more likely. You, however, will not
believe me. How shall I make my point? 

Let us assume I can measure the
strength and direction of the wind before every coin throw. I shall throw the
coin in the air many times, betting on heads or tails according to what I
reckon to be more likely, given the wind. If, at the end, I have collected
statistically significant winnings, you will have to acknowledge that I am in
the right.

Our situation is analogous. Instead of wind, we have $a_p$; instead of a coin,
we have whether or not $f(x)$ is a square mod $p$ for a random prime $x$.
(The prime $x$ stays fixed as $p$ varies.) If $f(x)$ (and thus $d$) lands on the more likely
side of squareness or non-squareness for significantly more than one-half of
all primes $p$, then $d$ will be sufficiently strange.

We can let $X_p$ be a random variable taking the value 
$\frac{-a_p}{p}$ when $f(x)$ is a
square mod $p$, and the value $\frac{a_p}{p}$ when $f(x)$ is a non-square mod
$p$. Then $X_p$ is $\frac{-a_p}{p}$ with probability $\frac{1}{2} -
\frac{a_p}{2 p}$, and $\frac{a_p}{p}$ with probability
$\frac{1}{2} + \frac{a_p}{2 p}$. It can be seen easily that
the expected value of $\sum_{p\leq z} X_p$ is
$\sum_{p\leq z} a_p^2/p^2 \sim \log \log z$.
We may assume pairwise independence and obtain that $\Var(\sum_{p\leq z} X_p)
= \sum_{p\leq z} \left(\frac{a_p^2}{p^2} - \frac{a_p^4}{p^4}\right)
 \sim \log \log z - O(1) \sim \log \log z$.
Thus, Chebyshev gives us that
$\sum_{p\leq z} X_p$ is $> (1- o(1)) \log \log z$ a
proportion
$1$ of the time. Now let $Y_p$ be $\frac{-a_p}{p}$ 
when a random integer $d\leq x$ is
a square residue modulo $p$; let $Y_p$ be $\frac{a_p}{p}$ otherwise. What is the probability that 
$\sum_{p\leq z} Y_p$ be larger than $(1 - o(1)) \log \log z$?

The probability that $Y_p$ take either of its two possible values is $1/2$.
Suppose that the variables $Y_p$ were mutually independent. Then the
expected value $\mathbb{E}(e^{\sum Y_p})$ 
of $e^{\sum Y_p}$ would be the product of the expected values of
$e^{Y_p}$. We can use this as follows. First of all,
\[
\Prob(\sum_{p\leq z} Y_p > (1 - o(1)) \log \log z) \leq
\Prob(e^{\sum_{p\leq z} Y_p} > (\log z)^{1 - o(1)}) \leq
\frac{\mathbb{E}(e^{\sum_{p\leq z} Y_p})}{(\log z)^{1 - o(1)}}.\]
Now, as we were saying,
\[\begin{aligned}
\mathbb{E}(e^{\sum_{p\leq z} Y_p}) &= 
\mathbb{E}(\prod_{p\leq z} e^{Y_p}) = 
\prod_{p\leq z} \mathbb{E}(e^{Y_p}) = 
\prod_{p\leq z} \left(\frac{1}{2} e^{\frac{a_p}{p}} + \frac{1}{2}
e^{-\frac{a_p}{p}}\right)\\
&= \prod_{p\leq z} \left(1 + \frac{1}{2} \frac{a_p^2}{p^2} +
\frac{1}{4!} \frac{a_p^2}{p^4} + \dotsc\right) \ll
\prod_{p\leq z} e^{\frac{1}{2} \frac{a_p^2}{p^2}}\\
&=e^{\sum_{p\leq z} \frac{a_p^2}{p^2}}
= e^{(1 + o(1)) \log \log z} = (\log z)^{1 + o(1)} .
\end{aligned}
\]
Hence
\begin{equation}\label{eq:homero}
\Prob(\sum_{p\leq z} Y_p > (1 - o(1)) \log \log z) 
\ll 
\frac{1}{(\log z)^{1/2 - o(1)}} ,\end{equation}
which is the bound we desire.

Now, the variables $Y_p$ are not in fact mutually independent, and, since the
probabilities we are dealing with are close to $1/2$ rather than to $0$, we
cannot apply the tricks in Erd\"os-Kac. A simpler approach will in fact do.
Let $z = N^{\frac{1}{3 \log \log N}}$ and $k = \frac{1}{2} \log \log z$.
Then, while the variables $Y_p$ are not mutually independent, they are more
than pairwise independent: any $2k$ of them are mutually independent (with a
small error term). We can thus proceed as in the proof of Chebyshev's theorem,
taking a $(2k)$th power instead of a square. The bound thus obtained is
essentially as good as (\ref{eq:homero}): we obtain that
the probability that
$\sum_{p\leq z} Y_p/\sqrt{p}$ be larger than $(1 - o(1)) \log \log z$ 
is $O(1/(\log z)^{1/2 - o(1)}) = O(1/(\log N)^{1/2 - o(1)})$.

We conclude that, if $d$ is strange in the two ways
we have considered  -- having
numbers of prime divisors differing from the norm, and ``agreeing with
the wind'' for considerably more than half of all $p$'s -- then it lies
in a set of cardinality at most
\[O\left(N \cdot \frac{1}{(\log N)^{\frac{1}{2} \log 3 + \frac{1}{2} -
      \epsilon}}\right) = O\left(\frac{N}{(\log N)^{1.0493\dotsc- \epsilon}}\right) .\]
This is smaller than $o\left(\frac{N}{\log N}\right)$,
which was the goal we set ourselves in the discussion after (\ref{eq:malato}).
Thus (\ref{eq:malato}) is indeed at most $o(N/\log N)$, and we are done
proving the main theorem.

\end{document}